# About norms, semi-norms and variation of trigonometric splines

Denysiuk VP Dr of phys-math. sciences, Professor, Kiev, Ukraine

National Aviation University

kvomden@nau.edu.ua

## Annotation

The work examines norms in $L^2$ of fundamental trigonometric splines of odd and even degrees, which in some cases coincide with polynomial ones. Fundamental trigonometric splines for the case where the convergence factors depend on the parameter are also considered. It is shown that an important characteristic of these splines of odd and even degrees is their complete variation. The semi-norms of trigonometric splines of odd degrees are also studied. The given material is illustrated by graphs; so, in particular, an illustration of the property of the least curvature of polynomial splines is given for the first time.

## Introduction

Approximation, respectively representation, of an arbitrary known or unknown function through a set of some special functions can be considered as a central topic of analysis. We will use the term "special functions" to refer to classes of algebraic and trigonometric polynomials and their modifications; at the same time, we believe that the classes of trigonometric polynomials also include trigonometric series. As a rule, such special functions are easy to calculate and have interesting analytical properties [1].

One of the most successful modifications of algebraic polynomials is polynomial splines, which are stitched together from segments of these polynomials according to a certain scheme. The theory of polynomial splines appeared relatively recently and is well developed (see, for example, [2], [3], [4], etc.). The advantages of polynomial splines include their approximate properties [5]. The main disadvantage of polynomial splines is their lumpy structure, which greatly complicates their use in analytical transformations.

Later it turned out[6], [7] that there are also modifications of trigonometric series that depend on several parameters and have the same properties as polynomial splines [8] [9]; moreover, the class of such modified series is quite broad and includes the class of polynomial periodic splines. This gave reason to call the class of such series trigonometric splines.

Convergence of trigonometric series that define trigonometric splines, is provided by convergence factors [6], [7], which have a decreasing order $O\left(k^{-(1+r)}\right)$, ($r = 0, 1, ...$); because at $r \geq 1$ these series coincide uniformly, then they are trigonometric Fourier series of a special form.

Important characteristics of trigonometric of splines are their norms in the functional space $L^2_{[0,2\pi]}$ [5] and half-norms [2]. Also, in our opinion, it is appropriate to introduce such a characteristic of these splines as their variation into consideration. It is clear that all these characteristics of the norm are largely determined by the type of selected convergence factors and the vector $\Gamma$, on which trigonometric splines depend. Therefore, it is relevant to study these characteristics of trigonometric splines with different convergence factors and vectors $\Gamma$.

## The purpose of the work

Comparison of standards, variations and half-norms of trigonometric splines of arbitrary degrees with some convergence factors and vectors $\Gamma$.

## The main part

Trigonometric splines with Riemann convergence factors were considered in [9,12]. $\sigma(r,k)$, power factors of convergence $v(r,k)=(k)^{-(1+r)}$, and some other convergence factors, when applying which we obtain different trigonometric splines; these splines interpolate a given sequence of values of some periodic function at the nodes of interpolation grids. Also, these splines at the same parameter values $r$ ($r = 0,1,2,...$) have the same differential properties and belong to space $2\pi$-periodic, continuous functions having continuous derivatives up to order $r-1$ including; we denote such a space by $C^{r-1}$. Note that the symbol $C^{-1}$ we denote the set $2\pi$-periodic, piecewise constant functions with a finite number of discontinuities of the first kind.

Note that all splines with the listed convergence factors at $r \to \infty$ converge to interpolation trigonometric polynomials that interpolate a given sequence of values of some periodic function at the nodes of interpolation grids.

In the further study of such trigonometric splines, it is important to consider their norms in space $L^2$, which are defined as follows.

$$\|f\| = \left( \int_0^{2\pi} [f(t)]^2 \, dt \right)^{1/2}. \quad (1)$$

Note that such norms of simple polynomial splines were studied in detail in [5],

When comparing the norms of trigonometric splines of type (1), we chose to compare the norms of fundamental trigonometric splines [11]. Recall that the fundamental spline on a uniform grid $\Delta_N^{(I)} = \{t_j^{(I)}\}_{j=1}^N$, where $I$ - grid indicator, ($I = 0,1$), $t_j^{(0)} = \frac{2\pi}{N}(j-1)$, $t_j^{(1)} = \frac{\pi}{N}(2j-1)$, is called a spline $s_k^{(I)}(t)$, which is determined from the condition

$$s_k^{(I)}(t_j^{(I)}) = \begin{cases} 1, & k = j; \\ 0, & k \neq j, \end{cases} \quad k, j = 1, 2, ..., N.$$

Using fundamental splines $s_k^{(I)}(t)$, ($k = 1, 2, ..., N$), interpolation spline $S^{(I)}(f,t)$, which interpolates the function $f(t)$ in grid nodes $\Delta_N^{(I)}$, can be written in the form

$$S^{(I)}(f,t) = \sum_{k=1}^N f(t_k^{(I)}) s_k^{(I)}(t).$$

The choice of fundamental splines for the comparison of their norms is explained by the fact that, firstly, these norms of these splines do not depend on the values of the interpolated function, and secondly, such splines suggest a visual comparison of their shape, which, in our opinion, is important in applications.

In [10], a general method for constructing trigonometric fundamental splines of a degree was given $r$, ($r = 1, 2, ...$), with a stitching grid $I_1$, interpolation grid $I_2$, a vector parameters $\Gamma$ and power-sign constant convergence factors $\sigma 0(\alpha, r, k) = \alpha(k)^{-(1+r)}$, ($k = 1, 2, ...$). This method is defined by formulas

$$st_k^{(I_1, I_2)}(\Gamma, \sigma 0, \alpha, r, N, t) = \frac{1}{N}\left( 1 + 2 \sum_{j=1}^{\frac{N-1}{2}} \frac{c_j^{(I_1, I_2)}(\Gamma, \sigma 0, \alpha, r, k, N, t)}{h_j^{(I_1, I_2)}(\Gamma, \sigma 0, \alpha, r, k, N)} \right), \quad (2)$$

where

$$c_j^{(I_1, I_2)}(\Gamma, \sigma 0, \alpha, r, k, t) = \gamma_1 \sigma 0(\alpha, r, j) \cos(j(t - (x_k^{(0)})^{1-I_2}(x_k^{(1)})^{I_2})) +$$

$$+ \sum_{m=1}^{\infty} (-1)^{m(I_1+I_2)} \left[ \gamma_2 \sigma 0(\alpha, r, mN - j)(-1)^{1+r} \cos((mN-j)(t - (x_k^{(0)})^{1-I_2}(x_k^{(1)})^{I_2})) + \quad (3)\right.$$

$$\left. + \gamma_3 \sigma 0(\alpha, r, mN + j) \cos((mN+j)(t - (x_k^{(0)})^{1-I_2}(x_k^{(1)})^{I_2})) \right],$$

$$h_j^{(I_1, I_2)}(\Gamma, \sigma 0, \alpha, r, N) = \gamma_1 \sigma 0(\alpha, r, j) + \sum_{m=1}^{\infty} (-1)^{m(I_1+I_2)} \left[ \gamma_2 (-1)^{1+r} \sigma 0(\alpha, r, mN - j) + \gamma_3 \sigma 0(\alpha, r, mN + j) \right]. \quad (4)$$

Note that the convergence factor $\sigma 0(\alpha, r, k)$, ($k = 1, 2, ...$), is the simplest sign-constant factor of convergence with decreasing order $O(k^{-(1+r)})$, and the trigonometric spline with such a convergence factor has degree $r$ and belongs to space $C^r$.

So, we will consider trigonometric fundamental splines (2) in the case when $\Gamma = \{\gamma_1, \gamma_2, \gamma_3\}$; while parameters $\gamma_k$ ($k = 1, 2, 3$) take arbitrary real values and $\gamma_2$ and $\gamma_3$ do not rotate simultaneously by 0. This requirement is explained by the fact that in the case when $\gamma_1 = 1$, a $\gamma_2 = \gamma_3 = 0$, we have the usual fundamental trigonometric polynomials [11].

If the vector $\Gamma$ looks like $\Gamma = \{1, 1, 1\}$, then we called fundamental trigonometric splines simple trigonometric fundamental splines. At the same time, functions $st_k^{(I_1, I_2)}(\Gamma, \sigma, \alpha, r, N, t)$ $h_j^{(I_1, I_2)}(\Gamma, \sigma, \alpha, r, N, k)$ we will mark accordingly $st_k^{(I_1, I_2)}(\sigma, \alpha, r, N, t)$ and $h_j^{(I_1, I_2)}(\sigma, \alpha, r, N, k)$.

When selecting a vector $\Gamma = \{\gamma_1, \gamma_2, \gamma_3\}$ it is advisable to take such considerations into account. It is clear that expression (3) can be conditionally divided into three components: low-frequency, medium-frequency and high-frequency; while the component $\gamma_1$ determines the influence of the low-frequency component on the trigonometric spline, and the components $\gamma_2$ i $\gamma_3$ determine the influence of the medium-frequency and high-frequency components of this spline, respectively. It is clear that it is possible to formulate a number of problems regarding the selection of vector components $\Gamma = \{\gamma_1, \gamma_2, \gamma_3\}$, based on the conditions of minimization of certain functional; we will not delve into this problem here.

In this work, we will choose a vector for illustration $\Gamma$ so that its components differ significantly from the components of a simple trigonometric spline; at also, we will consider the case when the vector $\Gamma$ looks like $\Gamma = \{.1, .5, 1.5\}$; splines with such a vector retain their notation $st_k^{(I_1, I_2)}(\Gamma, \sigma, \alpha, r, N, t)$.

It should be noted that the norms of fundamental trigonometric splines $st_k^{(I_1, I_2)}(\sigma, \alpha, r, N, t)$ and $st_k^{(I_1, I_2)}(\Gamma, \sigma, \alpha, r, N, t)$ do not depend on the index $k$, ($k = 1, 2, ..., (N-1)/2$); therefore, in the notation of the norms of these splines, we will omit this index. Also spline norms $st_k^{(0,0)}(\sigma, \alpha, r, N, t)$ and $st_k^{(0,1)}(\sigma, \alpha, r, N, t)$ comply with the norms $st_k^{(1,1)}(\sigma, \alpha, r, N, t)$ and $st_k^{(1,0)}(\sigma, \alpha, r, N, t)$; this also applies to splines $st_k^{(0,0)}(\Gamma, \sigma, \alpha, r, N, t)$ and $st_k^{(0,1)}(\Gamma, \sigma, \alpha, r, N, t)$. Taking this into account, in the following we will give only the norms of splines $st_k^{(0,0)}(\sigma, \alpha, r, N, t)$, $st_k^{(0,1)}(\sigma, \alpha, r, N, t)$ and $st_k^{(0,0)}(\Gamma, \sigma, \alpha, r, N, t)$ and $st_k^{(0,1)}(\Gamma, \sigma, \alpha, r, N, t)$.

Let's move on to a direct consideration of the norms of simple fundamental splines for some convergence factors and different powers $r$, putting for certainty $N = 9$; recall that the parameter $N$ determines the number of mesh nodes $\Delta_N^{(I)}$.

1. Consider the convergence factor $\sigma 0(\alpha, r, k) = \alpha k^{-(1+r)}$ [10]. First of all, we note that in the process of calculations it was found that the value of the parameter $\alpha$ ($0 < \alpha < 2$) do not affect the norm of trigonometric fundamental splines. Taking into account this fact, in the notation of splines in the future we will omit the dependence of these splines on the parameter $\alpha$.

Table 1 shows the values of the norms of simple trigonometric fundamental splines $st^{(0,0)}(\sigma 0, r, 7, t)$ and $st^{(0,1)}(\sigma 0, r, 7, t)$; note that splines $st^{(0,0)}(\sigma 0, r, 7, t)$ have polynomial counterparts.

Table 1.

| Grid indices | Value of parameter r (degree of spline) | | | | | | | | |
|---|---|---|---|---|---|---|---|---|---|
| | 1 | 2 | 3 | 4 | 5 | 6 | 7 | 8 | 50 |
| 0.0 | .6029 | 1.5558 | .7852 | 1.1108 | .8292 | .9895 | .8548 | .9427 | .8976 |
| 0.1 | 2.6886 | .7750 | 1.2644 | .8143 | 1.0352 | .8436 | .9614 | .8646 | .8976 |

Table 2 shows the values of the norms of trigonometric splines $st^{(0,0)}(\Gamma, \sigma 0, r, 7, t)$ and $st^{(0,1)}(\Gamma, \sigma 0, r, 7, t)$.

Table 2.

| Grid indices | Value of parameter r (degree of spline) | | | | | | | | |
|---|---|---|---|---|---|---|---|---|---|
| | 1 | 2 | 3 | 4 | 5 | 6 | 7 | 8 | 50 |
| 0.0 | .3405 | 3.7126 | .6884 | 26.3147 | .8524 | 3.9411 | .77816 | 1.3913 | .8976 |
| 0.1 | 8.8973 | 0.8041 | 2.7841 | .7553 | 45.2270 | .7711 | 1.9323 | .7952 | .8976 |

The results shown in Table 1 require discussion.

First of all, we note that earlier we showed that at $r \to \infty$ simple fundamental interpolating trigonometric splines converge to an interpolating trigonometric polynomial for any decreasing convergence factors $(1+r)$. Therefore, the values of the norms given in the last column can be considered asymptotically the limiting values of the norms at $r \to \infty$.

On the grids of the same name (0, 0) (or (1, 1)) norms of simple fundamental trigonometric (therefore polynomial) splines of odd degree with growth $r$ go, monotonically increasing, to the limit value from below, and splines of even power, on the contrary, go to the limit value, monotonically decreasing, from above. Norms of splines behave in almost the same way $st^{(0,0)}(\Gamma, \sigma 0, r, 7, t)$, however, the decay of even-degree spline norms is not monotonic.

The situation is the opposite on variable grids (0,1) (or (1,0)); the norms of simple fundamental trigonometric (hence polynomial) splines of even degrees go towards the limit value, monotonically increasing, from below, and splines of odd degrees go towards this value, monotonically decreasing, from above.

Almost the same situation is observed with trigonometric norms splines $st^{(0,0)}(\Gamma, \sigma 0, r, 7, t)$ and $st^{(0,1)}(\Gamma, \sigma 0, r, 7, t)$, however, as before, there is no monotony of decreasing norms of splines of odd powers.

The sharp increase in trigonometric norms is somewhat incomprehensible splines $st^{(0,0)}(\Gamma, \sigma 0, 4, 7, t)$ and $st^{(0,1)}(\Gamma, \sigma 0, 5, 7, t)$; we attribute this to vector exposure $\Gamma$.

2. Let us now consider the norms of trigonometric fundamental splines $st^{(0,0)}(\sigma, r, N, t)$ $st^{(0,1)}(\sigma, r, N, t)$ and $st^{(0,0)}(\Gamma, \sigma, r, N, t)$ and $st^{(0,1)}(\Gamma, \sigma, r, N, t)$ with sign-changing convergence factors, the simplest of which is the factor

$$\sigma(\alpha, r, k) = [\operatorname{sinc}(\alpha k)]^{1+r} = \left[\frac{\sin(\alpha k)}{\alpha k}\right]^{1+r}.$$

Trigonometric fundamental splines and their order derivatives $q$, ($q = 0, 1, ..., r-1$), with such convergence factors are constructed according to the formulas

$$st_k^{(I_1, I_2)}(\Gamma, \sigma, \alpha, r, q, N, t) = \frac{1}{N}\left(I(q) + 2\sum_{j=1}^{\frac{N-1}{2}} \frac{c_j^{(I_1, I_2)}(\Gamma, \sigma, \alpha, r, q, k, N, t)}{h_j^{(I_1, I_2)}(\Gamma, \sigma, \alpha, r, k, N)}\right), \quad (5)$$

where

$$I(q) = \begin{cases} 1 & \text{if } q = 0 \\ 0 & \text{otherwise} \end{cases}$$

$$c_j^{(I_1, I_2)}(\Gamma, \sigma, \alpha, r, q, k, t) = \gamma_1 \sigma(\alpha, r, j) j^q \cos(j(t - (x_k^{(0)})^{1-I_2} (x_k^{(1)})^{I_2}) + q\frac{\pi}{2}) +$$

$$+ \sum_{m=1}^{\infty} (-1)^{m(r+1+I_1+I_2)} \left[ \gamma_2 \sigma(\alpha, r, mN-j)(mN-j)^q \cos((mN-j)(t - (x_k^{(0)})^{1-I_2} (x_k^{(1)})^{I_2}) + q\frac{\pi}{2}) + (6) \right.$$

$$\left. + \gamma_3 \sigma(\alpha, r, mN+j)(mN+j)^q \cos((mN+j)(t - (x_k^{(0)})^{1-I_2} (x_k^{(1)})^{I_2}) + q\frac{\pi}{2}) \right],$$

$$h_j^{(I_1, I_2)}(\Gamma, \sigma, \alpha, r, N) = \gamma_1 \sigma(\alpha, r, j) + \sum_{m=1}^{\infty} (-1)^{m(r+1+I_1+I_2)} \left[ \gamma_2 \sigma(\alpha, r, mN-j) + \gamma_3 \sigma(\alpha, r, mN+j) \right]. \quad (7)$$

These are the norms of splines depend on grid indices $I_1, I_2$, vector $\Gamma$, parameter $\alpha$, which is included in the convergence factor $\sigma(\alpha, r, k)$, parameter $r$, which determines the degree of the spline and the index $k$ ($k = 1, 2, ..., N$) of the grid node. First of all, we note that the norms of splines do not depend on the grid node index $k$. When considering the norms of trigonometric fundamental splines, we will limit ourselves to considering splines of both even and odd degrees. Note that splines $st^{(0,0)}(\sigma, \alpha, r, 0, 7, t)$ of odd degrees have polynomial analogy at $\alpha = \pi/N$; therefore, the results of Table 1 are applicable to them.

We will begin the further explanation by considering the graphs of trigonometric fundamental splines $st^{(0,0)}(\sigma, \alpha, r, 0, 7, t)$ and $st^{(0,1)}(\Gamma, \sigma, r, 0, N, t)$ for $r = 3, 5, 7$ and $r = 2, 4, 6$ at $\alpha = \pi/N$; note that for spline reduction $st^{(0, I_2)}(\sigma, r, 0, N, t)$ on the graphs are denoted as $st(0, I_2, \sigma, r, 0, N, t)$, and splines $st^{(0, I_2)}(\Gamma, \sigma, r, 0, N, t)$ - as $st1(0, I_2, \sigma, r, 0, N, t)$, where $I_2 = 0, 1$.

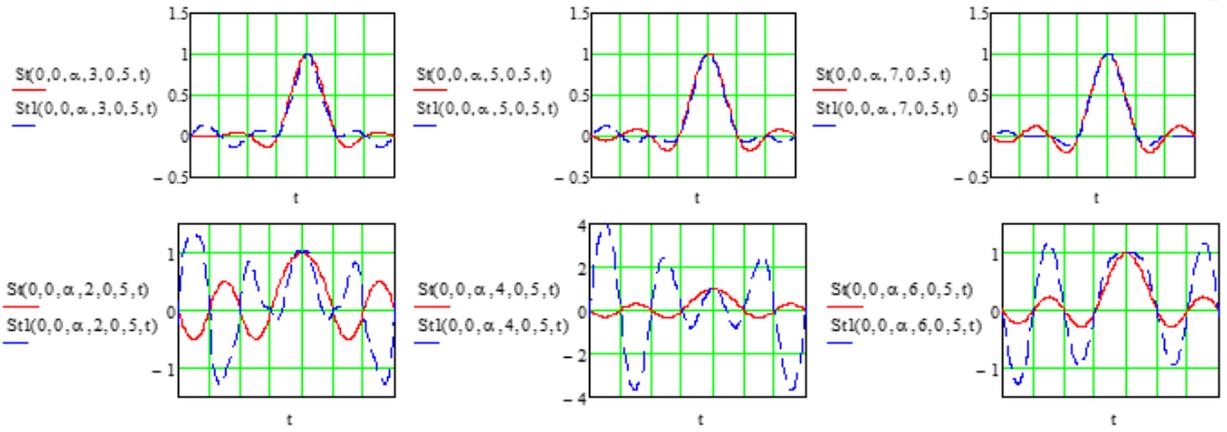

Fig. 1. Graphs of trigonometric splines $st^{(0,0)}(\sigma,\alpha,r,0,7,t)$ and $st^{(0,0)}(\Gamma,\sigma,\alpha,r,0,7,t)$ for $r=3,5,7$ and $r=2,4,6$ at $\alpha=\pi/N$.

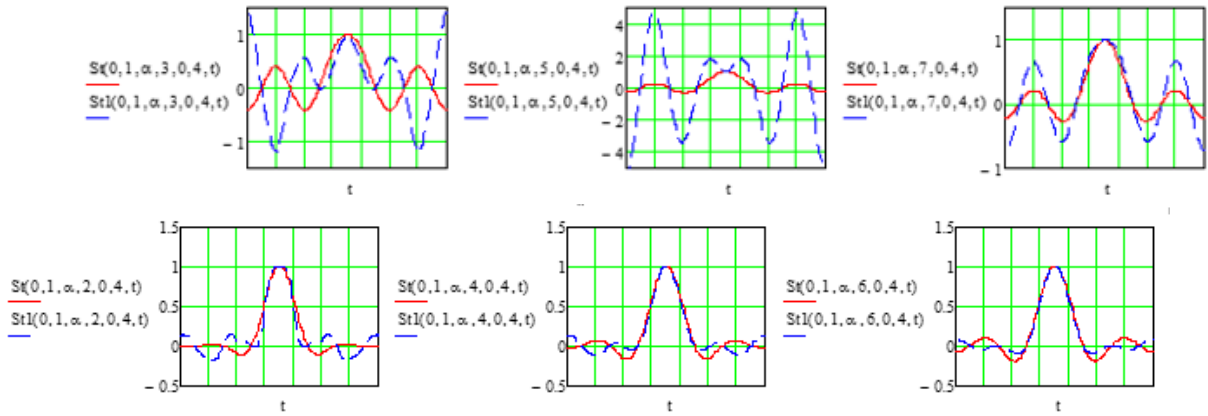

Fig. 2. Graphs of trigonometric splines $st^{(0,1)}(\sigma,\alpha,r,0,7,t)$ and $st^{(0,1)}(\Gamma,\sigma,\alpha,r,0,7,t)$ for $r=3,5,7$ and $r=2,4,6$ at $\alpha=\pi/N$.

Given the presence of significant spline fluctuations $st^{(0,0)}(\sigma,\alpha,r,0,7,t)$ and $st^{(0,0)}(\Gamma,\sigma,\alpha,r,0,7,t)$ even powers and splines $st^{(0,1)}(\sigma,\alpha,r,0,7,t)$ and $st^{(0,1)}(\Gamma,\sigma,\alpha,r,0,7,t)$ of odd degrees, in the future we will consider these splines only of odd and even degrees, respectively. Before plotting the graphs of the norms of these splines as functions of the parameter , it is necessary to make the following remark.

Graphs of norms and other characteristics of trigonometric splines as functions of a parameter α we will quote for the case $0<\alpha<\pi/2$; in addition, on these graphs, the value of these characteristics of the trigonometric spline, which has a polynomial counterpart, is indicated by a dashed line. The graphs of the splines themselves are given for $t\in[0,2\pi]$.

Let us present graphs of the norms of these splines $st^{(0,0)}(\sigma,\alpha,r,0,7,t)$ and $st^{(0,1)}(\Gamma,\sigma,r,0,N,t)$ as parameter functions α for $r=3,5,7$ and $r=2,4,6$; we will denote these norms accordingly $P0I_2(\alpha,r)$ and $P10I_2(\alpha,r)$.

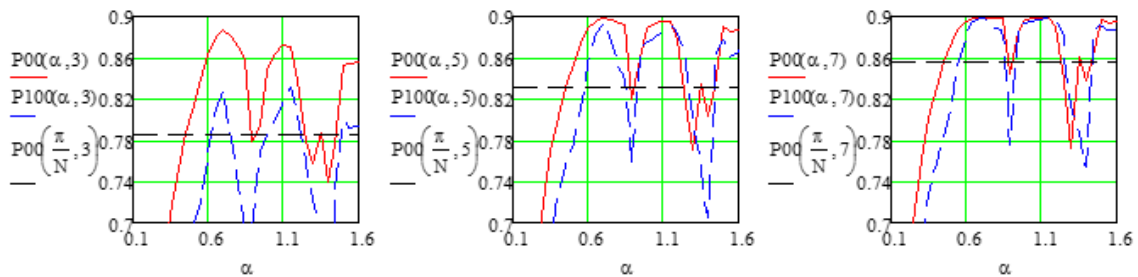

Fig. 3. Graphs of norms of trigonometric splines $st^{(0,0)}(\sigma,\alpha,r,0,7,t)$ and $st^{(0,0)}(\Gamma,\sigma,\alpha,r,0,7,t)$ at $r=3,5,7$.

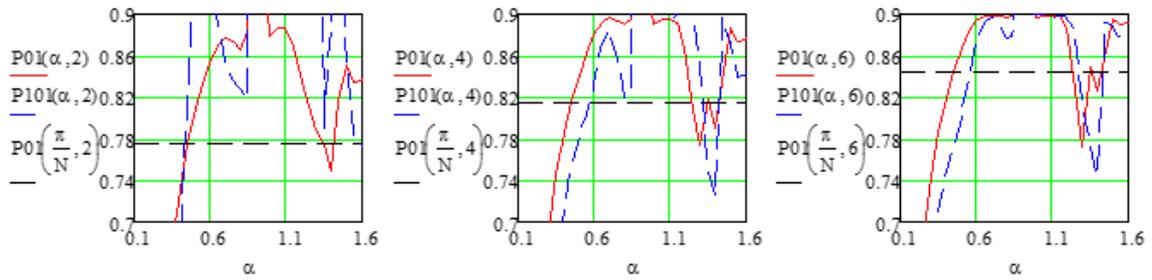

Fig. 4. Graphs of norms of trigonometric splines $st^{(0,1)}(\sigma,\alpha,r,0,7,t)$ and $st^{(0,1)}(\Gamma,\sigma,\alpha,r,0,7,t)$ at $r = 2, 4, 6$.

From the analysis of graphs of trigonometric norms of the splines shown in Figs. 3, 4, it follows that there are such values of the parameter $\alpha$, in which the norms of trigonometric splines are smaller than the norms of polynomial splines. Consider the graphs of such trigonometric splines, for example, at $\alpha = .15$.

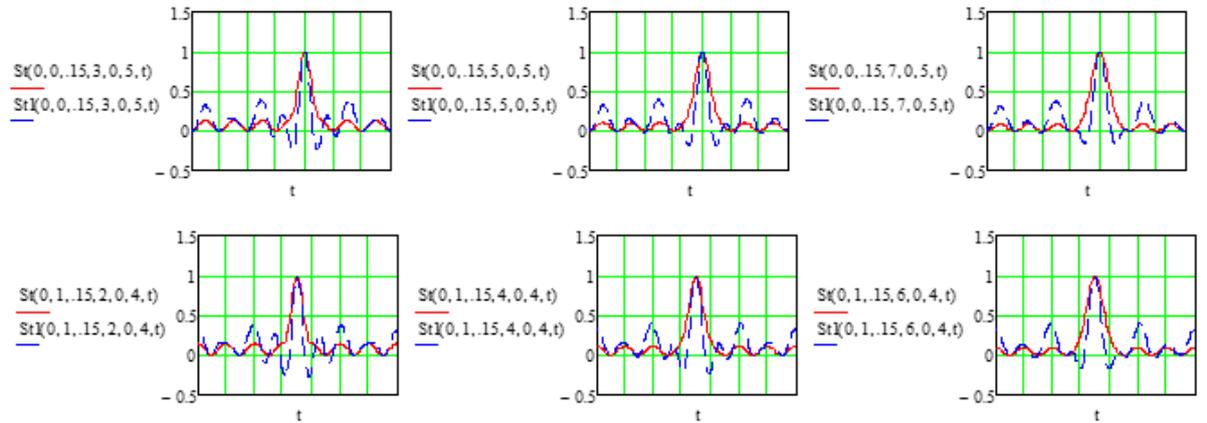

Fig. 5. Graphs of trigonometric splines $st^{(0,I_2)}(\sigma,\alpha,r,0,7,t)$ and $st^{(0,I_2)}(\Gamma,\sigma,\alpha,r,0,7,t)$ at $\alpha = .15$, $r = 3,5,7$ and $r = 2,4,6$; $I_2 = 0,1$.

Comparing the forms of trigonometric graphs splines shown in Fig. 1, 2, obtained at $\alpha = \pi/N$ with the forms of the graphs of these splines shown in Fig. 5, obtained at $\alpha = .15$, it is easy to come to the conclusion that at $\alpha = .15$ undesirable fluctuations of trigonometric splines increase significantly. Thus, it can be concluded that the norms of trigonometric splines do not fully characterize these splines; it is necessary to introduce some other characteristic of trigonometric splines that would numerically characterize the shape of these splines.

The shape of trigonometric splines can be numerically characterized in at least two ways.

With one of them, the presence of fluctuations splines can be characterized by the variation of these splines over the interval of their change. This approach can be applied to trigonometric splines of both even and odd degrees.

The second approach is research seminorms of trigonometric splines and application of consequences from Holladay's theorem (Holladay JC) [2]. A certain disadvantage of this approach is that it is applicable only for splines of odd powers. Let's consider both approaches in more detail.

As you know, the variation of a function on a segment $[a,b]$ is defined as follows. $\Delta = \{x_i\}_{i=1}^{N}$. $a = x_0 < x_1 < ... < x_N = b$, - any partition of the segment $[a,b]$. In aria (also a complete variation or complete change) of a function $f(x)$ on a segment $[a,b]$ is called a quantity

$$V_a^b(f) = \sup_\Delta \sum_{i=1}^{N} |f(x_{i-1}) - f(x_i)|. \quad (8)$$

It should be noted that for functions with a continuous first derivative (namely, such functions include trigonometric splines at $r \geq 2$) variation can be calculated using the formula

$$V_a^b(f) = \int_a^b |f'(x)| dx. \quad (9)$$

In addition, for such functions, the length of the arc of splines, which is calculated by the formula, can be considered

$$L(f) = \int_a^b \sqrt{1+(f'(x))^2}\,dx. \quad (10)$$

It is clear that formulas (9) and (10) lead to different results; but we are not interested in the absolute values of the variation or arc length, but only in the value of the parameter $\alpha$, at which certain values of variation or arc length are reached.

We will present graphs of the variation of trigonometric splines of even and odd powers; with $V_0^{2\pi}\left[st^{(0,I_2)}(\sigma,\alpha,r,0,7,t)\right]$ we will mark $V0I_2(\alpha,r)$, a $V_0^{2\pi}\left[st^{(0,I_2)}(\Gamma,\sigma,\alpha,r,0,7,t)\right]$ - by $V10I_2(\alpha,r)$.

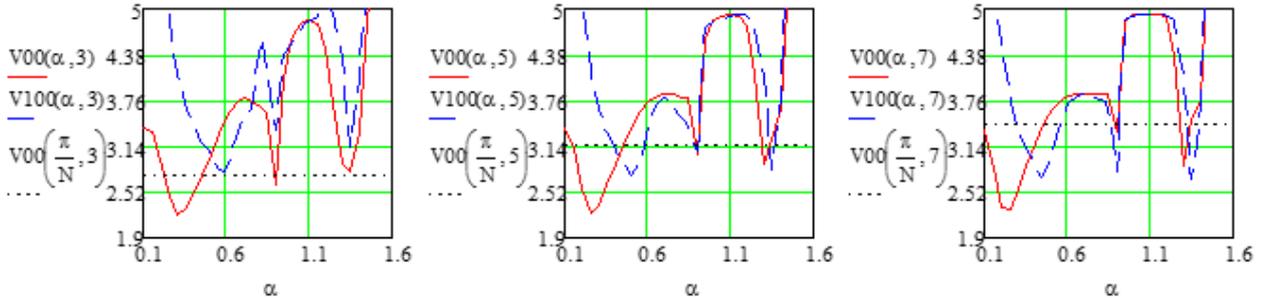

Fig. 6. Variation graphs of trigonometric splines $st^{(0,0)}(\sigma,\alpha,r,0,7,t)$ and $st^{(0,I_2)}(\Gamma,\sigma,\alpha,r,0,7,t)$ at $r=3,5,7$.

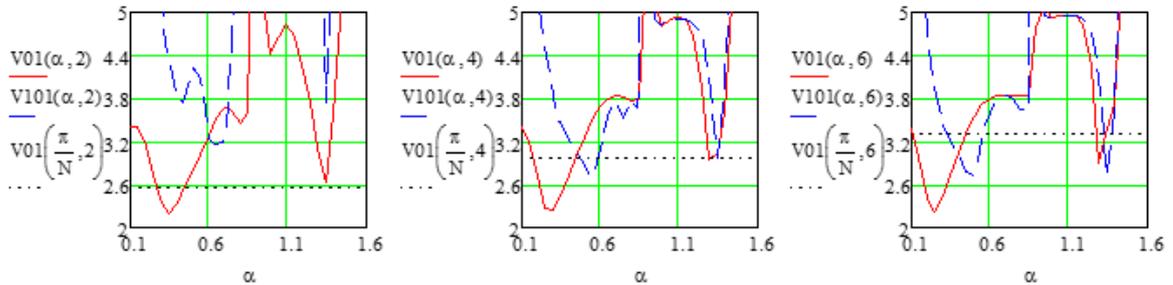

Fig. 7. Variation graphs of trigonometric splines $st^{(0,0)}(\sigma,\alpha,r,0,7,t)$ and $st^{(0,I_2)}(\Gamma,\sigma,\alpha,r,0,7,t)$ at $r=2,4,6$.

We will also give trigonometric onessplines with parameter values $\alpha$, at which the smallest variation is achieved (these values are approximate and were selected from the graphs in Fig. 6, 7); these splines are compared with trigonometric splines of the same powers that have polynomial counterparts.

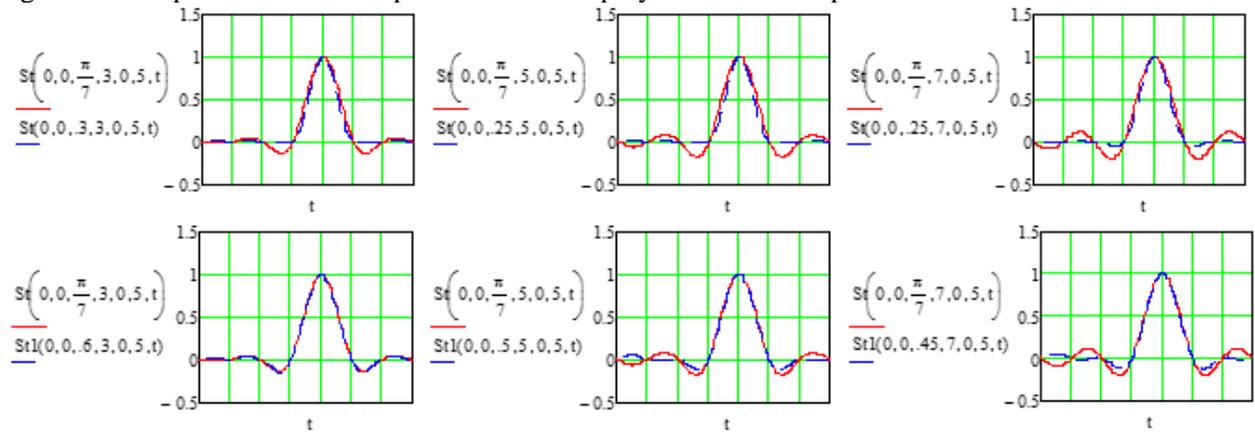

Fig. 8. Graphs of trigonometric splines $st^{(0,0)}(\sigma,\alpha,r,0,7,t)$ and $st^{(0,0)}(\Gamma,\sigma,\alpha,r,0,7,t)$ at $r=3,5,7$.

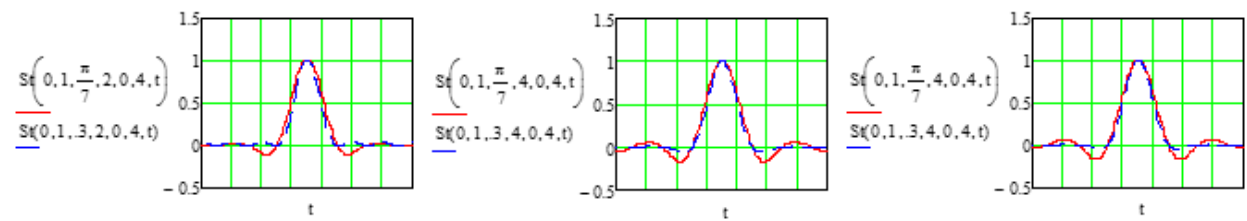

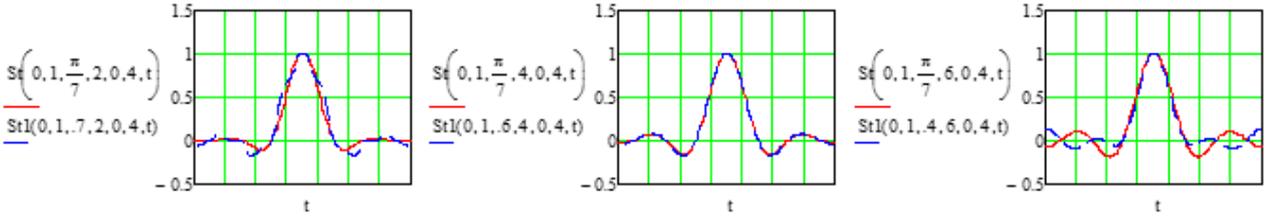

Fig. 9. Graphs of trigonometric splines $st^{(0,0)}(\sigma,\alpha,r,0,7,t)$ and $st^{(0,I_2)}(\Gamma,\sigma,\alpha,r,0,7,t)$ at $r = 2, 4, 6$.

It follows from the above graphs that only in spline $st^{(0,I_2)}(\Gamma,\sigma,.65,2,0,7,t)$ the variation exceeds the variation of the trigonometric spline $st^{(0,0)}(\sigma,\pi/7,2,0,7,t)$, which has a polynomial counterpart; in all other cases, the trigonometric splines have less variation than the corresponding trigonometric splines $st^{(0,I_2)}(\sigma,\pi/7,2,0,7,t)$, which, as we have already said, have polynomial counterparts. In addition, it turned out that the trigonometric spline $st^{(0,0)}(\Gamma,\sigma,.\alpha,3,0,7,t)$ and $st^{(0,1)}(\Gamma,\sigma,.\alpha,5,0,7,t)$ due to the selection of the parameter $\alpha$ can be made quite close to trigonometric splines $st^{(0,0)}(\sigma,\pi/7,3,0,7,t)$ and $st^{(0,1)}(\sigma,\pi/7,5,0,7,t)$ in accordance.

Let us now consider the second approach. As is known [2], polynomial simple interpolation splines of odd degree $2k-1$ ($k = 1, 2, ...$) have the property of minimal curvature; this property is that of all the class functions $C^k$ only on these splines the minimum value of the functional (half-norm) is reached, which has the form

$$\|f\|_k = \left(\int_0^{2\pi} \left[f^{(k)}(t)\right]^2 dt\right)^{1/2}.$$

We present the graphs of half-norms of splines of odd degree for different values of the parameters $\alpha$ and $r$; at the same time half-norms $st^{(0,0)}(\sigma,\alpha,r,0,7,t)$ and $st^{(0,0)}(\Gamma,\sigma,\alpha,r,0,7,t)$ we will mark accordingly $PP00(\alpha,r)$ and $PP100(\alpha,r)$.

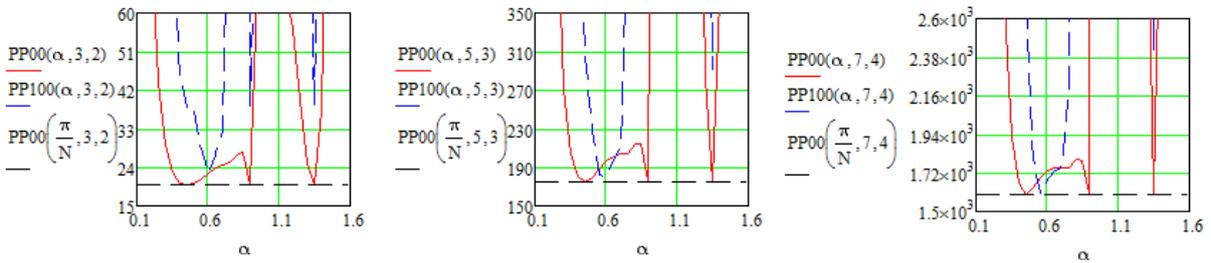

Fig. 10. Graphs of semi-norms of trigonometric splines $st^{(0,0)}(\sigma,\alpha,r,0,7,t)$ and $st^{(0,0)}(\Gamma,\sigma,\alpha,r,0,7,t)$ at $r = 3, 5, 7$.

Analysing the graphs shown in Fig. 10, the following conclusions can be drawn. First of all, let's note that we got a great illustration properties of least curvature of polynomial splines $st^{(0,0)}(\sigma,\alpha,r,0,7,t)$ odd degree, since they coincide with trigonometric splines of the same degree at $\alpha = \pi/N$.

It also turned out that trigonometric splines have the property of least curvature (and therefore coincide with polynomial splines) and when $\alpha = k\pi/N$, ($k$ and $N$ - mutually prime numbers; $k \neq N$, $k = \pm 1, \pm 2, ...$). In other words, there is a countable set of parameter values $\alpha$, in which trigonometric splines $st^{(0,0)}(\sigma,\alpha,r,0,7,t)$ coincide with simple polynomial splines of odd degree. Trigonometric splines $st^{(0,0)}(\Gamma,\sigma,\alpha,r,0,7,t)$ do not have the property of least curvature, and therefore not at any values of the parameter $\alpha$ do not coincide with polynomials. We also note that there are the following parameter values $\alpha$ in which trigonometric splines $st^{(0,0)}(\sigma,\alpha,5,0,7,t)$, $st^{(0,0)}(\sigma,\alpha,7,0,7,t)$ and $st^{(0,0)}(\Gamma,\sigma,\alpha,5,0,7,t)$ and $st^{(0,0)}(\Gamma,\sigma,\alpha,7,0,7,t)$ have the same values of semi-norms. There are also such values of this parameter for which the semi-norms of splines $st^{(0,0)}(\Gamma,\sigma,\alpha,5,0,7,t)$ and $st^{(0,0)}(\Gamma,\sigma,\alpha,7,0,7,t)$ are smaller than the half-norms of the splines $st^{(0,0)}(\sigma,\alpha,5,0,7,t)$, $st^{(0,0)}(\sigma,\alpha,7,0,7,t)$.

Finally, we note that for trigonometric splines $st^{(0,0)}(\sigma,\alpha,r,0,7,t)$ of even power there is also a countable set of parameter values $\alpha$, in which they coincide; this set has the form $\alpha = (2k+1)\pi/N$, where $2k+1$ and $N$ - mutually prime numbers, ($k = 0, \pm 1, ...$).

## Conclusions

1. The paper examines the classes of fundamental trigonometric splines, which do not depend on the approximated function and also assume a visual comparison of their shape.
2. The norms of fundamental trigonometric splines are considered $st^{(0,I_2)}(\sigma 0, r, N, t)$ and $st^{(0,I_2)}(\Gamma, \sigma 0, r, N, t)$ of even and odd degree with sign constant convergence factors $\sigma 0(r)$ for $I_2 = 0, 1$.
3. It is shown that on the grids of the same name (0, 0) (or (1, 1)) the norms of trigonometric simple (hence also polynomial) splines $st^{(0,0)}(\sigma 0, r, 7, t)$ odd power with growth $r$ go monotonically increasing to the limiting value from below, and splines of even degree, on the contrary, go monotonically decreasing to the limiting value from above. The same situation is observed with spline norms $st^{(0,0)}(\Gamma, \sigma 0, r, 7, t)$, however, the decay of even-degree spline norms is not monotonic.
4. The situation is the opposite on variable grids (0,1) (or (1,0)); the norms of trigonometric simple (hence, polynomial) splines of even degrees go towards the limit value monotonically increasing from below, and splines of odd degrees go towards this value monotonically decreasing from above. The same situation is observed with the norms of trigonometric splines $st^{(0,0)}(\Gamma, \sigma 0, r, 7, t)$ and $st^{(0,1)}(\Gamma, \sigma 0, r, 7, t)$, however, as before, there is no monotony of decreasing norms of splines of odd powers.
5. Dependencies of norms, semi-norms and variations of trigonometric fundamental splines are studied $st^{(0,I_2)}(\sigma, \alpha, r, q, N, t)$ and $st^{(0,I_2)}(\Gamma, \sigma, \alpha, r, q, N, t)$ of even and odd power with sign-changing convergence factors $\sigma(\alpha, r)$ for $I_2 = 0, 1$ from the parameter $\alpha$, ($0 < \alpha < \pi/2$) and different vectors $\Gamma$.
6. It is shown that there is a countable set of parameter values $\alpha$ species $\alpha = k\pi/N$, ($|k|$ and $N$ - mutually prime numbers; $|k| \neq N$, $k = \pm 1, \pm 2, ...$), where are simple polynomial splines $st^{(0,0)}(\sigma, \alpha, r, 0, N, t)$ of odd degree coincide with polynomial ones.
7. It is shown that there is a countable set of parameter values $\alpha$ species $\alpha = (2k+1)\pi/N$, ($|2k+1|$ and $N$ - mutually prime numbers; $|2k+1| \neq N$, $k = 0, \pm 1, ...$), where are simple polynomial splines $st^{(0,1)}(\sigma, \alpha, r, 0, N, t)$ of even power coincide with polynomial ones.
8. It is shown that only the norms of trigonometric splines incompletely characterize splines $st^{(0,I_2)}(\sigma, \alpha, r, q, N, t)$ and $st^{(0,I_2)}(\Gamma, \sigma, \alpha, r, q, N, t)$.
9. For a given vector $\Gamma$, different grid types and parameter values $r$ parameter values exist $\alpha$, at which the variation takes the smallest values.
10. An illustration of the least curvature property of simple polynomial splines is obtained, since trigonometric splines $st^{(0,0)}(\sigma, \alpha, r, 0, 7, t)$ at the parameter values $\alpha = k\pi/N$, ($k$ and $N$ - mutually prime numbers; $k \neq N$, $k = 1, 2, ...$), coincide with polynomials.
11. For a given vector $\Gamma$, different types of grids and odd values of the parameter $r$ values exist $\alpha$, at which the half-norms of trigonometric splines take the smallest values; these values $\alpha$ differ from the values at which the smallest variation is achieved.
12. In our opinion, based on research results, it is advisable to introduce classes of trigonometric splines of odd and even degrees with the smallest variation for a given vector $\Gamma$ and different types of grids; the application of splines of this class in approximation theory problems requires further research.
13. Also, in our opinion, it is advisable to introduce classes of trigonometric splines of odd degree with the smallest curvature for a given vector $\Gamma$ and different types of grids.
14. Undoubtedly, the issues discussed in the work require further research.

# List of references